\documentclass[11pt]{article}
\usepackage{amsmath,amsfonts,amssymb,amsthm,amscd}
\def\title#1{\def\thetitle{#1}}
\def\author#1{\def\theauthor{#1}}
\def\date#1{\def\thedate{#1}}

\title{Cohomology of Courant algebroids with split base}
\author{Gr\'egory Ginot\footnote{\texttt{ginot@math.jussieu.fr}\\ UPMC Paris 6\\
  Equipe Analyse Alg\'ebrique \\ 4, place Jussieu 75252 Paris, France} and 
  Melchior Gr\"utzmann\footnote{\texttt{grutzman@math.psu.edu}\\
     109 McAllister Building\\
     University Park, PA 16802\\USA\\}
}

\date{May 13, 2008}

\newcommand{\noprint}[1]{}  
\def \MG{\noprint}
\def \MGI{\noprint}

\def \GGi{\noprint}
\def \PX{\noprint}
\long\def\notneeded{\noprint}  
\long\def\deprecated{\noprint}

\newtheorem{vdef}{Definition}[section]
\newtheorem{lem}[vdef]{Lemma}
\newtheorem{prop}[vdef]{Proposition}
\newtheorem{cor}[vdef]{Corollary}
\newtheorem{thm}[vdef]{Theorem}

\theoremstyle{definition}
\newtheorem{rem}[vdef]{Remark} 
\newtheorem{ex}[vdef]{Example}

\def\[{\begin{equation}}  
\def\]{\end{euqation}}
\DeclareMathOperator\Ad{Ad}
\DeclareMathOperator\Ann{Ann}
\def\B{\mathcal{B}}
\def\Lie{{\mathcal{L}}}
\def\C{{\mathcal{C}}}

\def\ud{{\mathrm{d}}}

\def\E{{\mathcal{E}}}
\def\g{{\mathfrak{g}}}
\DeclareMathOperator\im{im}
\def\img{\im}
\def\Lie{{\mathcal{L}}}
\def\pfrac#1#2{\frac{\partial#1}{\partial#2}}
\def\pt{{\mathrm{pt}}}

\def\R{{\mathbb{R}}}

\def\smooth{\ensuremath{\C^\infty}}
\def\gsmooth{\smooth}
\def\T{{\mathbb{T}}}
\def\X{{\mathcal{X}}}
\def\Xx{{\mathfrak{X}}}
\def\Z{{\mathbb{Z}}}
\def \Hal{H_{Lie}}  

\begin{document}
\begin{center}
\Large \textbf{\thetitle}\\
\large \theauthor\\
\thedate
\end{center}

\begin{abstract} In this paper we study  the cohomology $H^\bullet_{st}(E)$ of a Courant algebroid $E$. We prove that if $E$ is regular, $H^\bullet_{st}(E)$ coincides with the naive cohomology $H_{naive}^\bullet(E)$ of $E$ as conjectured by Sti\'enon and Xu~\cite{SX07}. For general Courant algebroids $E$ we define a spectral sequence converging to $H^\bullet_{st}(E)$. If $E$ is with split base,  we prove that there exists a natural transgression homomorphism $T_3$ (with image in $H^3_{naive}(E)$) which, together with $H_{naive}^\bullet(E)$, gives all $H^\bullet_{st}(E)$. For generalized exact Courant algebroids, we give an explicit formula for $T_3$ depending only on the {\v S}evera characteristic clas of $E$.  
\end{abstract}

\section{Introduction}
The purpose of this paper is to study the cohomology of Courant algebroids.
The Courant bracket was first introduced by T.~Courant in 1990
(see~\cite{Cour90}) in
order to describe Dirac manifolds, a generalization of presymplectic and
Poisson manifolds.  In 1997  Liu, Xu and Weinstein introduced the notion of a
Courant algebroid in order to describe Manin triples for Lie bialgebroids
(\cite{Xu97}). Recently, Courant algebroids have been used as  a background to describe generalized
complex geometry, see \cite{Hit03,Gua04} and  as target spaces for
$3$-dimensional topological field theory~\cite{Ikeda01,Ikeda02,Par00,
  HoPa04,Royt06}. 

Roughly speaking, a Courant algebroid is a pseudo-Euclidean vector bundle
$E\to M$ together with an anchor map $\rho:E\to TM$ and a bracket $[.,.]$ on
$\Gamma E$ which satisfy the basic identities, \emph{e.g.} skew-symmetry, Jacobi
identity, Leibniz rule and ad-invariance,  only up to anomalies (which are
\emph{exact} terms). Up to the anomalies, the bracket and the anchor map are
similar to those of a Lie algebroid. Indeed, Courant algebroids appear to be
the right framework for pseudo-metric vector bundles equipped with something
like an ad-invariant Lie algebra structure. A primary example of a Courant
algebroid 
is given by the double $A\oplus A^*$ of a Lie bialgebroid $A$ or a Lie
quasi-bialgebroid, that is a Lie bialgebroid twisted by a
$3$-form~\cite{Royt99}.  

An important step forward  by Roytenberg  was a description of  Courant algebroids in
terms of a derived
bracket as introduced by Kosmann-Schwarzbach in \cite{Kos96}, see
\cite{Royt02}, or equivalently in terms of a nilpotent odd 
operator (also known as $Q$-structure). Hence the Courant algebroid structure 
with its intricate axioms can all be encoded
in a cubic function $H$ on a graded symplectic manifold and its derived
bracket. To do so, one goes into the 
context of graded manifolds and considers the graded manifold $E[1]$. The
pseudo-metric on $E$ makes $E[1]$ into a (graded) Poisson manifold. By 
constructions of Weinstein, {\v S}evera~\cite{SevPC} and Roytenberg~\cite{Royt02}, there is a minimal symplectic realization $(\E,
\{.,.\})$ of $E[1]$. Now, there is a cubic Hamiltonian $H$
satisfying $\{H,H\}=0$ on $\E$ encoding the Courant algebroid structure
together with the symplectic structure of $\E$. For instance, the Courant
bracket is given by the formula $[\phi,\psi]=\{\{H,\phi\},\psi\}$.

The derived bracket construction also leads to a natural notion of cohomology
of a Courant algebroid. Since $\{H,H\}=0$, the operator
$Q=\{H,.\}:\gsmooth(\E)\to\gsmooth(\E)$ is a differential. Hence one can
define the cohomology of $E$ as the cohomology $H^\bullet(\gsmooth(\E),Q)$.
So far, there are only few examples of Courant algebroids for which the
cohomology is known. For instance when $E\cong T^*M\oplus TM$ is an exact
Courant algebroid, its cohomology is isomorphic to the de Rham cohomology of
$M$. On the other hand, if the base is a point, $E$ is a Lie algebra
(together with an ad-invariant pseudo-metric) and its cohomology is
isomorphic to its cohomology as a Lie algebra. One reason which makes
$H^\bullet(\gsmooth(\E),Q)$ rather difficult to treat is that its
construction relies on the minimal symplectic realization $\E$ and not just
on $E$ or $E^*$ itself. In particular it is quite different from the usual
cohomology theories for ``Lie theoretic objects'' such as Lie algebroids or
Leibniz algebras where the cohomology is defined using a differential given
by a Cartan type formula. For instance, the (de Rham) cohomology of a Lie
algebroid $(A,[.,.],\rho:A\to M)$ over $M$ is the cohomology of the complex
of forms $(\Gamma(\Lambda^\bullet A^*),\ud_A)$, where the differential $\ud_A$
is given, for $\alpha$ a $n$-form and $\psi_1,\dots,\psi_n \in \Gamma(A)$, by
the Cartan formula
\begin{align} \label{eq:dintro} \langle\ud_A\alpha, \psi_1\wedge\ldots\psi_{n+1}\rangle :=
  &\sum_{i=1}^{n+1} (-1)^{i+1}\rho(\psi_i)\langle\alpha,\psi_1\wedge\ldots
    \widehat{\psi_i}\ldots\wedge\psi_n\rangle  \\
  &+\sum_{i<j} (-1)^{i+j} \langle\alpha,[\psi_i,\psi_j]\wedge\psi_1\wedge
    \ldots\widehat{\psi_i}\ldots\widehat{\psi_j}\ldots\psi_n\rangle  .\nonumber
\end{align}
This formula does not make sense for a Courant algebroid $E$ in place of $A$
since,  there is no definition of a Courant algebroid for which the bracket
can be made skew-symmetric and the Jacobi identity and the Leibniz rule are
strictly satisfied.  Nevertheless,  Sti\'enon-Xu~\cite{SX07} recently
observed that formula~\eqref{eq:dintro} makes sense  if one restricts $E$ to
the kernel $\ker \rho$ of its anchor map. Precisely, they proved that 
Formula~\eqref{eq:dintro} defines a differential $\ud_{naive}$ on $\Gamma(\Lambda^\bullet 
\ker \rho)$, where the pairing $\langle .,.\rangle$ in the formula is induced 
by the identification of $E$ and $E^*$ by the pseudo-metric. They call the 
cohomology of $(\Gamma(\Lambda^\bullet \ker \rho),\ud_{naive})$ the naive 
cohomology of $E$ and denoted it $H_{naive}^\bullet(E)$.  Furthermore, they 
conjectured that for a transitive Courant algebroid $E$, \emph{i.e.} a Courant 
algebroid with surjective anchor map, the naive cohomology of $E$ is 
isomorphic to the (standard) cohomology $H^\bullet(\gsmooth(\E),Q)$ of $E$. We 
prove this conjecture in Section~\ref{s:split}, see Corollary~\ref{c:trans}. A 
nice feature of the naive cohomology is that it can be calculated using the 
same techniques as for Lie algebroids since it is defined similarly.

Another goal of this paper is to study the cohomology of general Courant 
algebroids. We show that in general, the naive cohomology and the standard 
cohomology are related by a spectral sequence. Spectral sequences are a useful 
tool in geometry and topology. For instance, given a map $X\xrightarrow{\pi} 
Y$, the Leray spectral sequence allows one to compute the cohomology of $X$ in 
terms of the (sheaf) cohomology of $Y$ and the inverse images of $\pi$. This 
spectral sequence can be very complicated if $\pi$ is not nice enough. On the 
other hand, it takes a much simpler form, when $\pi$ is a fibration with a 
trivial action of $\pi_1(Y)$ on the fiber $F$, relating the cohomology of $X$ 
to the cohomology of $Y$ and $F$. Spectral sequence also leads to 
transgression~\cite{BT, SSeq85}. Our spectral sequence can be thought of as an 
analogue of the Leray spectral sequence of the morphism of graded manifolds 
$\E \to (\ker \rho)^*$ induced by the pseudo-metric and the map $\E \to E$. 
Indeed, the pseudo-metric allows to see sections of $\Gamma(\ker \rho)$, which 
are functions on $(\ker\rho)^*[1]$, as 
functions on $E[1]$ which can be pulled back to functions on $\E$. Then, the 
ideal $I=(\Gamma(\ker \rho))\subset \gsmooth(\E)$ is a differential ideal 
which induces a filtration of the complex $(\gsmooth(\E),Q)$. The filtration 
gives rise to a spectral sequence  converging to $H^\bullet(\gsmooth(\E),Q)
$ that we call the \emph{naive ideal spectral sequence}, see 
Section~\ref{s:naiveSS}.

\smallskip

The next important part of this paper is to compute explicitly the naive ideal
spectral sequence in terms of smooth geometric data. This computation
involves taking quotient of $TM$ by $D=\im \rho$, thus, one has to assume
this quotient is nice. In Section~\ref{s:split}, we consider the case of
Courant algebroids with split base, \emph{i.e.}, Courant algebroids for which
the base is isomorphic to a product $L\times N$ of manifolds and $D=\im
\rho\cong TL\times N$. \PX{makes no sense}\MGI{I think now it does}Note that any 
transitive Courant
algebroid is with split base. In the general case of split base, the naive 
cohomology is an
approximation of the cohomology of $E$ in the sense that the  term
 $E_2^{\bullet,\bullet}$ of the spectral sequence is the tensor product of
the naive cohomology with symmetric multivector fields on $N$, see
Proposition~\ref{p:E2}.

Furthermore, by transgression, the spectral sequence induces a map
$T_3:\X(N) \to H^3_{naive}(E)$ 
 that we call the \emph{transgression
  homomorphism}. This map is a cohomological characteristic class. Indeed, we
show that the naive cohomology and the transgression homomorphism determines
explicitly the (standard) cohomology of $E$, see Theorem~\ref{thm:split}.

There is a nice class of Courant algebroids parametrized by closed
$3$-forms. According to {\v S}evera~\cite{SevLett},  exact Courant algebroids
over $M$, \emph{i.e.}, Courant algebroids fitting in an exact sequence $0\to
T^*M \to E \to TM\to 0$ are in bijection with cohomology classes of closed
$3$-forms on $M$~\cite{SeWe01, Royt99}. We define a \emph{generalized exact
  Courant algebroid} to be a regular Courant algebroid fitting into an exact
sequence $0\to D^* \to E \to D\to 0$ of bundles over $M$, where $D=\im
\rho$. By {\v S}evera's argument, we get that generalized exact Courant
algebroids are in bijection with cohomology classes of closed $3$-forms on
(the Lie algebroid) $D$. The class associated to a generalized exact Courant
algebroid is called its {\v S}evera characteristic class. If, furthermore,
$E$ is with split base $D\cong TL\times N$, then a $3$-form on $D$ lies in
$\Omega^3(L) \otimes \smooth(N)$. We prove that for a generalized exact Courant
algebroid with split base, the transgression homomorphism is given by
the map $\X(N)\ni X\mapsto (1\otimes X)(C) \in H^3(L)\otimes
\smooth(N)\cong  H^3_{naive}(E)$ where $[C]$ is the {\v S}evera
characteristic class, see Proposition~\ref{P:D=C}. As a consequence,  we
obtain an explicit computation of the cohomology of $E$ in terms of $[C]$, see Corollary~\ref{c:genexact}.

The plan of the paper is as follows. In Section~\ref{s:defn}, we recall the
definition of Courant algebroids and their cohomology,
following~\cite{Royt99, Royt02}. In Section~\ref{s:naive}, we explain the
definition of the naive cohomology following~\cite{SX07}. We define the naive
ideal spectral sequence in Section~\ref{s:naiveSS}. Then we compute the
spectral sequence in the case of Courant algebroids with split base in
Section~\ref{s:split}. In particular we prove the Conjecture of
Sti\'enon-Xu, see Corollary~\ref{c:trans}.  We also define the transgression
homomorphism. The second main result in this section is
Theorem~\ref{thm:split} computing the cohomology of a Courant algebroid with
split base\notneeded{ in terms of its transgression homomorphism and naive
  cohomology}.
In section~\ref{s:gexsplit}, we study generalized exact Courant algebroids,
defining their {\v S}evera characteristic class and giving their
classification, Proposition~\ref{P:classification}. We then compute the
transgression homomorphism in terms of the characteristic class, see
Proposition~\ref{P:D=C}. For the reader's convenience we included an appendix
recalling the basic facts and definition of graded (and super) geometry that
we need. 

\subsection{Aknowledgement}
The  authors would like to thank Ping Xu for various comments on the earlier
versions of this paper.  \notneeded{Especially for pointing out possible problems with a 
formula like the one before equation~\eqref{eq:badquotient} in examples like
the irrational torus.}

\section{Courant algebroids and their cohomology}\label{s:defn}
In this section we recall the  notion of Courant algebroids and their
cohomology. We work in the context of graded manifolds as described in
Appendix~\ref{app:graded}. We refer to~\cite{Vor91, Vor01, Sev01b, Royt02} for more details.

The Courant bracket was originally defined as a skew-symmetric bracket
satisfying the Jacobi identity and Leibniz rule only up to some
anomalies~\cite{Cour90}. As was shown by Roytenberg in his
thesis~\cite{Royt99}, there is an alternative  equivalent formulation of
Courant algebroids, namely one with a bracket fulfilling a certain
form of the Jacobi identity (in fact the Leibniz rule for the bracket see
Axiom (1) in Definition~\ref{def:Courant}) but no longer skew-symmetric originally due to Dorfman~\cite{Dor87,Dor93}.   We
recall the latter definition since it is more convenient in our context and
also the more closely related to Loday algebroids as in~\cite{SX07}.
\begin{vdef}[Courant algebroid] \label{def:Courant} A Courant algebroid is a 
vector bundle $E\to M$  equipped with a symmetric non-degenerate bilinear 
product
  $\langle.,.\rangle: E\otimes_M E\to \R$, an $\R$-linear bracket\PX{Use
     $\circ$ instead of bracket, because it's more common and your
     ``bracket'' is not skewsymmetric}\MGI{I'd like to stick to this notation 
  $[.,.]$. My understanding is that bracket means biderivation, not just (graded)
  skewsymmetry.} on
  the sections of $E$ and a bundle map $\rho:E\to TM$, called the anchor map.
  These three operations have to satisfy to the
  following rules\,:
  \begin{align}
    [\phi,[\psi_1,\psi_2]] &= [[\phi,\psi_1],\psi_2]
       +[\psi_1,[\phi,\psi_2]]  \label{ax:Jacobi}\\
    [\phi,f\psi] &= \rho(\phi)f\cdot\psi +f[\phi,\psi]  \label{ax:Leibn}\\
    [\phi,\phi] &=\frac12\rho^*\ud\langle\phi,\phi\rangle  \label{ax:nSkew}\\
    \rho(\phi)\langle\psi,\psi\rangle &= 2\langle[\phi,\psi],
      \psi\rangle  \label{ax:adInvar}
  \end{align} where $\phi,\psi,\psi_1,\psi_2$ are sections into $E$
  and $f$ is a smooth function on $M$.
\end{vdef}
\begin{rem}
Note that the bracket $[.,.]$ is not skew-symmetric. Also, some authors denote this bracket $\circ$, for instance see~\cite{SX07, Royt02}.  
\end{rem}
\begin{rem}
Since the bilinear form $\langle.,.\rangle$ is non-degenerate, a Courant algebroid is a pseudo-Euclidean vector bundle. In particular, there is a canonical identification of $E$ and $E^*$ induced by the pseudo-metric.
\end{rem}

\notneeded{\begin{rem}\PX{Notion of Courant algebroid is standard in
  literature and you are not the first to explain that.}
Note that the first two axioms are similar to those  of a Lie algebroid.  The third axiom shows that the bracket is  not skew-symmetric. Precisely $[\phi,\psi]-[\psi,\phi]$ is equal to an ``exact'' term $\rho^*\ud\langle\phi,\psi\rangle$ where $\rho^*: T^*M\to E$ is the map  dual to the anchor map $\rho$ given by  the identification of $E$ and $E^*$ by the pseudo-metric $\langle.,.\rangle$. Note that the original definition of Courant algebroid~\cite{Cour90} can be recovered by skew-symmetrizing the bracket, see~\cite{Royt99}.
\end{rem}}

A Courant algebroid is called \emph{transitive} if the anchor map $\rho:E\to
TM$ is surjective.  It is called \emph{regular} if $\rho$ is of constant rank.

Due to an observation of K.~Uchino (\cite{Uchi02}), the Jacobi identity
(Definition~\ref{def:Courant}.\ref{ax:Jacobi}) and Leibniz rule
(Definition~\ref{def:Courant}.\ref{ax:Leibn}) imply that the anchor map
preserves the bracket:
for all $\psi_1, \psi_2 \in \Gamma(E)$, one has
$$  \rho[\psi_1,\psi_2]
  =[\rho\psi_1,\rho\psi_2]_{TM}  $$
Therefore, in the regular case the image of the anchor map is a Lie algebroid.
\medskip

The Courant bracket can also be obtained as a derived bracket~\cite{Kos96, Royt02}
on a (degree 2) graded symplectic manifold. More precisely, the pseudo-metric
on $E$ induces a structure of (degree 2) graded Poisson manifold on
$E[1]$. \notneeded{Note that assigning odd degree to the fibers indeed makes the Poisson
bracket symmetric on the (odd) fibers.}  There is a minimal symplectic
realization $\E$ of $(E[1],\langle.,.\rangle)$ constructed as follows. When
$E$ is a split vector bundle $E=A\oplus A^*$ (with the obvious pairing), $\E$
is just $T^*[2]A[1]$ with the usual graded symplectic structure (see
Appendix~\ref{app:graded}). In the general case (where $\langle.,.\rangle$
might not even have split signature), a construction of $\E$ was found by
Weinstein, {\v S}evera~\cite{SevPC} and Roytenberg~\cite{Royt02}. First consider the graded
symplectic manifold $T^*[2]E[1]$. \notneeded{It is of too high dimension to be the
expected symplectic realization of $E$.} Let $\imath:E[1]\to E\oplus E^*[1]$ be
the isometric embedding (given by the pseudo-metric). The minimal symplectic
realization of $E[1]$ is the pullback
$\E:=\imath^*T^*[2]E[1]=E[1]\times_{E\oplus E^*[1]} T^*[2]E[1]$ of
$T^*[2]E[1]$ along $\imath$. We denote $\pi:\E\to E[1]$ the canonical bundle
projection. Note that $\E$ fits into the following short exact sequence of
graded fiber-bundles over $M$\,:
\begin{equation}\label{eq:EE} 0\to T^*[2]M\to \E \xrightarrow{\pi} E[1]\to 0
\;.\end{equation}

Roytenberg~\cite{Royt02} proved that there is a cubic Hamiltonian $H$  on
$\E$ encoding the Courant algebroid structure on $E$. More precisely, the
Hamiltonian $H$ satisfies the following properties\,:
\begin{description}
\item[(1)] The Courant bracket is given by the derived bracket
  $[\phi,\psi]=\{\{H,\phi\},\psi\}$ where we identify sections $\phi,\psi$ of
  $E$ with fiber-linear functions on $E$ by the pseudo-metric
  $\langle.,.\rangle$.\footnote{The derived brackets go back to an idea of
  Kosmann-Schwarzbach~\cite{Kos96}}
\item[(2)] The anchor map is given, for any $\psi \in E$, by the formula  $\rho(\psi)=\{\{\psi,H\},.\}$.
\item[(3)] $H$ is nilpotent, \emph{i.e.},  $\{H,H\}=0$.
\end{description}
\begin{rem} Note that axioms (\ref{ax:Leibn}), (\ref{ax:nSkew}) and
  (\ref{ax:adInvar}) of a Courant algebroid (in Definition~\ref{def:Courant})
  now follow directly from the derived bracket construction. The first
  axiom, \emph{i.e.}, the Jacobi identity is equivalent to $\{H,H\}=0$.
\end{rem}

The derived bracket approach allows us to define a natural notion  of
cohomology for Courant algebroids.
The Hamiltonian $H$ gives a degree 1 derivation $Q=\{H,-\}$ which is of
square zero (since $\{H,H\}=0$). Furthermore, $Q$ maps graded functions on
$\E$ to graded functions on $\E$. We denote $A^\bullet:=C^\infty(\E)$ the
graded functions on the minimal symplectic realization $\E \stackrel{\pi}\to
E$ of $E$, where, by definition $A^n$ stands for the functions of degree
$n$. Then $(A^\bullet,Q)$ is a complex.

Therefore, the following definition due to Roytenberg~\cite{Royt02} makes
sense\notneeded{ (recall that $\E\stackrel{\pi}\to E$ is a minimal realization of
$E[1]$)}.

\begin{vdef}[Cohomology of Courant algebroids] \label{def:Q} Let $(E\to M, 
[.,.] , \rho ,$ $ \langle .,.\rangle)$ be a Courant algebroid. Let $H$ be a function 
of degree 3 on the symplectic realization $\E$ generating the Courant bracket. 
The cohomology of the Courant algebroid $E$ is the cohomology of the complex 
$(A^\bullet=C^\infty(\E),Q)$  equipped with the differential 
$Q:=\{H,.\}$.
\end{vdef}
We  denote by $H^\bullet_{std}(E)$ the above cohomology of the Courant algebroid $E$.  Note that
$(A^\bullet,Q)$ is a differential graded commutative algebra, thus
$H^\bullet_{std}(E)$ is a graded commutative algebra.

Let us further introduce coordinates on $\E$ for later use.  We denote $x^i$
the  coordinates on the base $M$, $p_i$ their conjugates of degree 2 and
$\xi^a$ the (pseudo) orthogonal fiber-coordinates of degree 1 on $E$.  Then
the cubic Hamiltonian  reads in coordinates as 
$$H=\rho^i_a(x)p_i\xi^a +\frac16C_{abc}(x)\xi^a\xi^b\xi^c \;,$$  where
$\rho(\xi_a)=\rho^i_a(x)\pfrac{}{x^i}$ encodes the anchor map and
$C_{abc}:=\langle[\xi_a,\xi_b],\xi_c\rangle$ are the structure functions of
the bracket.

\section{Naive cohomology}\label{s:naive}
In this section, we recall the definition of the naive cohomology of a
Courant algebroid~\cite{SX07}. It is  less involved than
Definition~\ref{def:Q}; for instance it does not use a symplectic realization
of $E$.

Mimicking the definition of the differential giving rise to the cohomology of Lie   algebroids, the idea is to consider an operator $\ud: \Gamma(\Lambda^\bullet E)\to \Gamma(\Lambda^{\bullet +1} E)$ given    by
the Cartan formula\,:
\begin{align} \label{eq:d} \langle\ud\alpha, \psi_1\wedge\ldots\psi_{n+1}\rangle :=
  &\sum_{i=1}^{n+1} (-1)^{i+1}\rho(\psi_i)\langle\alpha,\psi_1\wedge\ldots
    \widehat{\psi_i}\ldots\wedge\psi_{n+1}\rangle  \\
  &+\sum_{i<j} (-1)^{i+j} \langle\alpha,[\psi_i,\psi_j]\wedge\psi_1\wedge
    \ldots\widehat{\psi_i}\ldots\widehat{\psi_j}\ldots\psi_{n+1}\rangle  \nonumber
\end{align}
where $\psi_1,\dots \psi_{n+1}$ are sections of $E$ and $\alpha\in 
\Gamma(\Lambda^n E)$ is identified with an $n$-form on $E$, \emph{i.e.}, a 
section of $\Lambda^n E^*$ by the pseudo-metric. However, the 
formula~\eqref{eq:d} is not well defined because it is not $\C^\infty(M)
$-linear in the $\psi_i$ (due, for instance, to axioms~\ref{ax:Leibn} and 
\ref{ax:nSkew} in Definition~\ref{def:Courant}.
 The operator $\ud$ does not square to zero either since it is not skewsymmetric
in the $\psi_i$.\footnote{Using the skew-symmetric bracket does not help either, 
because it only fulfills a modified Jacobi identity.}  Nevertheless, 
Sti{\'e}non-Xu~\cite{SX07} noticed that
the formula~\eqref{eq:d} for $\ud$ becomes $\C^\infty(M)$-linear in the
$\psi_i$ when one restricts to $\alpha\in\Gamma(\Lambda^n(\ker\rho))$.\footnote{
In this case using the skewsymmetric or Jacobi-fulfilling bracket does not matter
since the difference is exact, thus vanishes in the inner product with $\ker\rho$.} 
In fact they proved the following
\begin{lem}\label{lem:d}
 Formula~\eqref{eq:d} yields a well defined operator $\ud:\Gamma(\Lambda^\bullet \ker\rho)\to \Gamma(\Lambda^{\bullet +1} \ker\rho)$. Moreover, one has $\ud \circ \ud =0$.
\end{lem}
Note that  $\ker \rho$ may be a singular vector bundle.\footnote{In this
 case $\Gamma(\ker \rho)$ means smooth
 sections into $E$ that are pointwise in the kernel of $\rho$.  We define
 similarly sections of  $\Lambda^\bullet\ker\rho$ where  $\rho$ has been
 extended as an odd $\smooth(M)$-linear derivation $\Lambda^\bullet E\to
 \Lambda^\bullet E\otimes TM$.}

\begin{proof} 
The first claim follows 
from the fact that the failure of the Leibniz rule in the left-hand side
is an exact term, \emph{i.e.}, is in the image of $\rho^*$ and that $\rho
 \circ \rho^*=0$.  Now the terms for $\ud\circ\ud$ add up to zero using
the Jacobi identity as in the Lie algebroid case, since all terms are equivalent
to those using the skew-symmetric bracket. See~\cite{SX07} Section 1 for more details.
\end{proof}

\begin{rem}\label{rem:deRham}
For a Lie algebroid $A$, the space $(\Gamma(\Lambda^\bullet A^*),\ud)$, where
$\ud$ is the operator given by formula~\eqref{eq:d}, defines its cohomology.  
In particular, it calculates the de Rham cohomology of $M$ when $A=TM$.
\end{rem}

Thanks to the pseudo-metric $\langle.,.\rangle $ on $E$, one
can view $\Gamma(\Lambda^n E) $ as graded functions on $E[1]$ (of degree n),
which can further be pulled back to the minimal symplectic
realization $\pi: \E\to E[1]$.  Thus, we can identify
$\Gamma(\Lambda^n(\ker\rho))$ with a subalgebra of $\smooth(\E)$.  Sti\'enon-Xu proved~\cite{SX07} the following Proposition.
\begin{prop}\label{lem:Q=d}
 The Q-structure $Q=\{H,.\}$ (see Definition~\ref{def:Q}) maps
 $\Gamma(\Lambda^\bullet(\ker\rho))$ to itself. Moreover, if $\alpha \in \Gamma(\Lambda^\bullet(\ker\rho))$, then $Q(\alpha)=\ud(\alpha)$.
\end{prop}
In other words, $Q$ restricted to $\Gamma(\Lambda^\bullet(\ker\rho))$ coincides
with the differential $\ud$ given by formula~\eqref{eq:d} (note that
Proposition~\ref{lem:Q=d} also implies that $\ud\circ \ud=0$).
Since $\ud$ squares to 0,  Sti{\'e}non-Xu defined\,:

\begin{vdef}[Naive cohomology]\label{d:naive}  Let $(E,[.,.],\rho,\langle.,.\rangle)$ be a
Courant algebroid.  The naive cohomology of $E$ is the cohomology of  the sections of
$\Lambda^\bullet\ker\rho$ equipped with the differential $\ud$ given by  the Cartan-formula~\eqref{eq:d}.
\end{vdef}
We denote $H_{naive}^\bullet(E)=H^\bullet(\Gamma(\Lambda^\bullet\ker\rho),\ud)$
the naive cohomology groups of $E$.  By Proposition~\ref{lem:Q=d}, there is a
canonical morphism $\phi:H_{naive}^\bullet(E)\to H^\bullet_{std}(E) $ from the
naive cohomology to the  (standard) cohomology of Courant algebroids, see~\cite{SX07}. We will
prove that this morphism is an isomorphism in the transitive case, see
Corollary~\ref{c:trans}.

\section{Geometric spectral sequence for cohomology of Courant algebroids with
 split base}\label{s:regcase}
In this section we define a spectral sequence converging to the cohomology of a
Courant algebroid.  Then, in the case of Courant algebroids with split base,
we compute the spectral sequence in terms of \notneeded{smooth} geometrical
data.  For details about spectral sequences refer to \cite{CE56, SSeq85}.

\subsection{The naive ideal spectral sequence}\label{s:naiveSS}
The algebra $A^\bullet:=C^\infty(\E)$ of graded functions on $\E$ is endowed
with a natural filtration induced by the ideal generated by the kernel of the
anchor map $\rho: E\to TM$.   More precisely, let $I$ be the ideal
$$I:=\Gamma(\Lambda^{\ge1}\ker\rho)\cdot C^\infty(\E)\;,$$
{\it i.e.},  the ideal of functions containing at least one coordinate of
$\ker\rho$, where we identify sections of $E$ to odd functions on $E$ by the
pseudo-metric (as in Section~\ref{s:naive}). Since $\ker \rho$ gives rise to
the naive cohomology, we call $I$ the naive ideal of $\E$.
\begin{lem} $I$ is a differential ideal of the differential graded algebra
$(A^\bullet,Q)$
\end{lem}
\begin{proof}
According to Proposition~\ref{lem:Q=d} and Lemma~\ref{lem:d}, we have  $Q(\Gamma(\Lambda^n \ker\rho)) \subset\Gamma(\Lambda^{n+1} \ker \rho)$. The result follows since $Q$
is a derivation.
\end{proof}

Since $I$ is a differential ideal, we have a decreasing bounded (since $E$ is finite dimensional) filtration of differential graded algebras
$A^\bullet= F^0A^\bullet \supset F^1A^\bullet \supset F^2A^\bullet \dots$, where
$F^pA^q:= I^p\cap A^q$. Therefore\,:
\begin{prop}\label{l:SS}
There is a  spectral sequence of algebras
$$E_0^{p,q}:=F^pA^{p+q}/F^{p+1}A^{p+q} \Longrightarrow H^{p+q}_{std}(E) $$
converging to the cohomology of the Courant algebroid  $E$.
\end{prop}
We call this spectral sequence, the \emph{naive ideal spectral sequence}.
\begin{proof} The spectral sequence is the one induced by the filtration $F^\bullet A^\bullet$ of the complex $(A^\bullet,Q)$. It is convergent because the filtration is bounded.
\end{proof}
\begin{rem}
In order for the naive ideal spectral sequence to be useful, one needs to be
able to calculate the higher sheets $E_k^{p,q}$ of the spectral sequence in
terms of (smooth) geometry of $E$. Such calculations involve the image of the
anchor map. Thus it seems reasonable to restrict to the class of regular
Courant algebroids. \GGi{For general regular Courant algebroids,
the foliation groupoid given by the integral leaves of $D$ may be complicated
and suggests to make the computation in the framework of stacks, see
Remark~\ref{rem:stack} below.}\MG{It might not even be a stack, because the
groupoid is not Lie either.  See my notes after Remark~\ref{rem:stack}.}
\end{rem}

\subsection{Courant algebroids with split base} \label{s:split}

In this section we compute explicitly the naive ideal spectral sequence
of Proposition~\ref{l:SS} for what we call Courant algebroid with split base
and then prove the conjecture of Sti\'enon-Xu as a special case.

\medskip

By the bracket preserving-property of the anchor-map $\rho$, $D:=\img\rho$ is an integrable
(possibly singular) distribution.
\begin{vdef}[split base]\label{d:split}  A Courant algebroid $(E\to M,\langle.,.\rangle,
 [.,.],\rho)$ is said to have split base iff $M\cong L\times N$ and the image of
 the anchor map is $D:=\im\rho\cong TL\times N\subset TM$.
\end{vdef}
In particular, a Courant algebroid with split base is a regular Courant
algebroid.  Furthermore, the integral leaves of the distribution $D$ are
smoothly parametrized by the points of $N$; thus the quotient $M/D$ by the
integral leaves of the foliation is isomorphic to $N$.

\begin{rem} \label{rem:lcoo}
The local coordinates $\xi^a$ for $E$ introduced in Section~\ref{s:defn} can 
be splitted accordingly to the isomorphism $D\cong TL\times N$. This splitting 
is useful in order to do local computations. More precisely, over a 
chart-neighborhood $U$ of $M$, $D:=\img\rho$ can be spanned by coordinate 
vector fields $\pfrac{}{x^I}$. Let $\xi^I:=\rho^*\ud x^I$ ($\in \ker \rho$) be 
vectors that span the image of $\rho^*$, and let $\xi_I$ be preimages of 
$\partial_I$, dual to the $\xi^I$. Then choose coordinates $\xi^A\in \ker 
\rho$ normal (with respect to the pseudo-metric) and orthogonal to both the 
$\xi^I$s and the $\xi_I$s. Therefore we have split the coordinates $\xi^a$s in 
the three subsets consisting of the $\xi^I$s, the $\xi_I$s and the $\xi^A$s. 
Furthermore this splitting also induces a splitting of the degree 2 
coordinates $p_i$s (the conjugates of the coordinates on $M$) into the $p_I$s, 
which are the symplectic duals of the $x^I$, and the $p_{I'}$s (their 
complements for which the Poisson bracket with the $x^I$s vanish). Since 
$D\cong TL\times N$, the coordinates $x^I$ and $x^{I'}$ can be chosen to be 
coordinates of $L$ and $N$ respectively.

The Hamiltonian in these coordinates reads as $H=p_I\xi^I+\frac16C_{abc}(x)
\xi^a\xi^b\xi^c$.
Note that in order to compute structure functions like $C_{abc}$ you need
$\xi_A$ which is the dual frame of $\xi^A$ or due to the pseudo orthonormality
$\xi_A=\pm\xi^A$.
\end{rem}

\medskip

We denote $\Xx^q(N)$ the space of (degree q) symmetric multivector fields 
$\Gamma_N(S^{q/2}(TN))$ with the convention that $S^{q/2}(TN)$ is $\{0\}$ for 
odd $q$'s, \emph{i.e.}, $\Xx^\bullet(N)$ is concentrated in even degrees. With 
these notations, the sheet $E_1^{\bullet,\bullet}$ of the spectral sequence of 
Lemma~\ref{l:SS} is given by\,:
\begin{lem}\label{l:E1} For a Courant algebroid with split base $D\cong TL\times N$,
  one has
\begin{align*}
  E_1^{p,q}\cong \Gamma_M\Lambda^p(\ker\rho)\otimes \Xx^q(N).
\end{align*}
\end{lem}

\begin{proof}
The differential $\ud_0:E_0^{p,q}\to E_0^{p,q+1}$  is the
differential induced on the associated graded $\bigoplus
 F^pA^{p+q}/F^{p+1}A^{p+q}$ of the naive filtration $F^\bullet
 A^\bullet$. Hence, it is obtained from $Q$ by neglecting all terms which
contain at least one term in $I$.  Using the local coordinates of
Remark~\ref{rem:lcoo}, we find $\ud_0 = p_I\pfrac{}{\xi_I}$.  Globally, the
algebra $E_0^{0,\bullet}$ is isomorphic to $\gsmooth(B)$ for the graded
manifold $B:=\E/I$ which, on a local chart $U\subset M$ is isomorphic to
$B_{|U}\cong T^*[2]U\otimes_U D[1]|U$.  Globally, $B$ fits into the short exact
sequence of graded fiber bundles over $M$\,:
\begin{equation}\label{eq:Bseq}  0\to T^*[2]M\to B\to D[1]\to 0\;,
\end{equation} where the latter map
  $B\to D[1]$ is the map induced by the anchor map $\rho: E \to D=\im\rho$
on the quotient of $\E$ by $I$.
The differential $\ud_0$ on $E_0^{0,\bullet}\cong \smooth(B)$ is canonically 
identified with the odd vector field $\tilde\rho_0$ induced by $Q$ on the 
quotient $B=\E/I$.

There is a similar interpretation of $E_0^{1,\bullet}$. Precisely, $E_0^{1,
\bullet}$ is isomorphic to $\Gamma_B(B_1)$ for a (graded) vector bundle 
$B^1\to B$ which, locally, is the fiber product $B^1_{|U}\cong 
\ker\rho[1]\times_U B$. In particular, $B^1$ fits into the short exact 
sequence of graded fiber bundles over $M$\,:
\begin{equation}\label{eq:B1seq} 0\to\ker\rho[1]\to B^1\to B\to 0
\end{equation}
\MG{This seemed surprising to me first, but the first map follows since
  $\Gamma(\ker\rho)$ span $I/I^2$ as module over $\gsmooth(B)$.}
The isomorphism $E_0^{1,\bullet}\cong \Gamma_B(B_1)$ identifies the 
differential $\ud_0:E_0^{1,\bullet}\to E_0^{1,\bullet+1}$ with the odd vector 
field $\tilde\rho: \Gamma_B(B^1)\to \Gamma_B(B^1)$ defined as the covariant 
derivative $\tilde\rho=\nabla_{\tilde\rho_0}$ along $\tilde{\rho}_0$ where $\nabla$ 
is a local connection on $B_1$ vanishing on a local frame $\xi^a$ of $B^1$. 
This is well defined, because the transition functions between such frames 
come from functions on $M$ and $\tilde{\rho_0}$ projected to $M$ vanishes (we 
extend to arbitrary sections of $B^1$ via Leibniz rule).

This identification of $E_0^{1,\bullet}$ extends to the other lines
$E_0^{p\geq 2,\bullet}$ of the spectral sequence easily.  Namely, there is an
isomorphism $E_0^{p,\bullet}\cong \Gamma_B(S^pB^1)$, where $S^p$ stands for
the graded symmetric (hence it is skew-symmetric since the fibers
$\ker\rho[1]$ are of odd degree) product over $B$. We extend $\tilde{\rho}$
to $ \Gamma_B(S^pB^1)$ by the Leibniz rule. Since, by Lemma~\ref{l:SS},
$\ud_0$ is a derivation, we also have the identification of $\ud_0:E_0^{p\geq
  2,\bullet}\to E_0^{p\geq 2,\bullet+1}$ and $\tilde{\rho}$.

\smallskip

According to Lemma~\ref{l:SS}, $\ud_0$ is a derivation, thus it is sufficient 
to compute the cohomology of the complex $(E_0^{0,\bullet},d_0)$.\MGI{It
  seems we had to compute $E^0$ and $E^1$} The 
sequence~\eqref{eq:Bseq} yields a morphism of sheaves $\eta: \smooth(B) \to 
\Gamma_M(S^\bullet(TM)[2]) \to \Gamma_M(S^\bullet(TM/D)[2])$. On a local chart,
 the complex $(E_0^{0,\bullet},d_0)$ is isomorphic to the Koszul complex of 
$\B^\bullet:=\Gamma_M(S^\bullet(TM)[2]\otimes\Lambda^\bullet(E/\ker\rho))$ 
with respect to the differential $\tilde\rho= p_I\pfrac{}{\xi_I}$ induced by 
the regular family given by the $p_I$s. Since the image of $\tilde \rho$ spans 
$D\cdot \B^\bullet$, the morphism of sheaves $\eta:\smooth(B) \to 
\Gamma_M(S^\bullet(TM/D)[2])$ is locally a quasi-isomorphism. Thus we have 
\begin{align*}H^q(E_0^{0,\bullet},d_0)&\cong \Gamma_M(S^\bullet(TM/D)[2]) \\ 
&\cong \Gamma_{N}(S^\bullet T(N)[2])\otimes_{\C^\infty(N)}\C^\infty(L\times N) 
\\ &\cong \X_{flat}^q(N) \otimes_{\C^\infty(N)}\C^\infty(M) \;,
\end{align*} where 
the second line follows from the fact that $E$ has split base $D\cong TL\times 
N$. 
\notneeded{ To see that this is also the global result, note that lines 
 2--3 are global identities. it remains to argue, that also the 1st line is 
 global.

 Clearly $H^2(E_0^{0,\bullet},\ud_0=\tilde\rho_0)$ are the sections into a 
 graded fiberbundle $C$ over $M$. Dividing out the image of $\tilde\rho_0$ 
 corresponds to taking the kernel of the right hand map in \eqref{eq:Bseq}, 
 which leads to the left hand term $T[2]M$. Restricting to the kernel of 
 $\tilde\rho_0$ corresponds to dividing out $D$. That is what the right hand 
 side of the first line says. 
 For the line $E_1^{1,\bullet}$ note that $\tilde\rho$ is a horizontal lift of
 $\tilde\rho_0$.  Therefore taking cohomology wrt.\ $\tilde\rho$ does not
 affect the fibers $\ker\rho[1]$ of $B^1\to B$, only the base $B$.  The change
 in the base is the computed result for $E_1^{0,\bullet}$.  In order to see
 that the fibers glue to the bundle $\ker\rho[1]$ over $M$ note that the
 transition functions in $B^1\to B$ actually come from $M$ (see \ref{eq:B1seq}).

 The other lines are generated by $E_0^{1,\bullet}$ so the result follows.
} The 
computation of the line $E_1^{1,\bullet}$ is similar since $\tilde\rho$ is a 
horizontal lift of $\tilde{\rho}_0$. In particular, $\tilde{\rho}$ does not 
act on the fibers $\ker\rho[1]$ of $B^1\to B$, but only on the base $B$. Hence 
$$H^q(E_0^{1,\bullet},d_0)\cong \Xx^q(N) \otimes \Gamma_M(\ker\rho).$$ Since 
the lines, $E_0^{p,0}$ are generated by (products of elements of) $E_0^{1,0}$ 
and $\ud_0$ is a derivation, the result follows.
\end{proof}

\begin{rem} \label{rem:stack}\PX{Shorten this remark\,!}
Lemma~\ref{l:E1} is the main reason to restrict to Courant algebroids with
split base.   In
general one can consider the quotient $M/D$  of $M$ by the integral leaves of the integrable distribution $D:=\img\rho$.  Let
$\C^\infty(M)^D$ be the space of smooth functions on $M$ constant
along the leaves
and $\X_{flat}(M/D)$  the space of derivations of $\C^\infty(M)^D$.  To describe the spectral sequence using 
smooth geometry,  we would like a formula of the form\,: $$E_1^{p,q}\cong 
\Gamma_M(\Lambda^p \ker \rho)\otimes S^{q/2}(\X_{flat}(M/D)). $$ A quick analysis of 
the proof of Lemma~\ref{l:E1} shows that this formula will hold if and only if 
we have the relation
\begin{equation}\label{eq:badquotient}\Gamma_M(TM/D)
  \cong \Gamma_{M/D}(T_{flat}(M/D))\otimes_{\C^\infty(M)^D}
\C^\infty(M) \;.
\end{equation}
 Courant algebroids with split base are a large class  for which relation~\eqref{eq:badquotient} holds.  However Relation~\eqref{eq:badquotient} does not hold for every regular
  Courant algebroid.
For instance, take the Lie algebroid $D$ underlying the irrational torus. That is  $M=\T^2$ is foliated by the action of
a non-compact one parameter subgroup of $\T^2$ and $D$ is the subbundle of
$TM$ associated to the foliation. The leaves are dense. Let $E=D\oplus D^*$ be  a generalized  exact Courant algebroid (as in Example~\ref{e:genexactsplitbase}). Note that $D$ is regular of rank $1$, thus $\Gamma_M(TM/D)$ is non zero but  $\smooth(M)^D\cong \R$, thus
$\X_{flat}(M/D)=0$.   Thus formula~\eqref{eq:badquotient} does not hold for  $E$.
\end{rem}

\bigskip

Recall from Section~\ref{s:naive} that the naive cohomology
$H_{naive}^\bullet(E)$ is the cohomology of the complex
$(\Gamma_M\Lambda^\bullet(\ker\rho),\ud)$.  The second sheet of the spectral sequence is
computed by $H_{naive}^\bullet(E)$\,:
\begin{prop} \label{p:E2} Let $E$ be a Courant algebroid with split base
  $D\cong TL\times N$. Then, one has an isomorphism of graded algebras
$$E_2^{p,q} \cong H_{naive}^p(E)\otimes \Xx^q(N). $$
\end{prop}
\begin{proof}
It is a standard fact of spectral sequences~\cite{CE56, SSeq85} that the
differential
$$\ud_1: E_1^{p,q}\cong H^q(F^pA/F^{p+1}A,\ud_0) \to E_1^{p+1,q}\cong
 H^q(F^{p+1}A/F^{p+2}A,\ud_0)$$
is the connecting homomorphism in the
cohomology long exact sequence induced by  the short exact sequence of
complexes $0\to F^{p+1}A/F^{p+2}A \to F^{p}A/F^{p+2}A \to F^{p}A/F^{p+1}A \to
0$.
 On a local chart, we obtain that $\ud_1$ is given by the following formula
\begin{align*}\ud_1
  =& \xi^I\pfrac{}{x^I}+C_{IA}^K\xi^I\xi^A\pfrac{}{\xi^K}
    +\frac12C_{AB}^K\xi^A\xi^B\pfrac{}{\xi^K}
    +\frac12C_{IJ}^A\xi^I\xi^J\pfrac{}{\xi^A} \\
  &+C_{AI}^B\xi^A\xi^I\pfrac{}{\xi^B}+\frac12C_{AB}^C\xi^A\xi^B\pfrac{}{\xi^C}
\end{align*}
where we use the local coordinates introduced in Remark~\ref{rem:lcoo}.
Now it follows from the isomorphism $E_1^{p,q}\cong
\Gamma_M\Lambda^p(\ker\rho)\otimes \Xx^q(N)$ given by Lemma~\ref{l:E1} and the
above formula for $\ud_1$ that  $$\ud_1=\ud\otimes 1:
\Gamma_M\Lambda^p(\ker\rho)\otimes \Xx^q(N) \to
\Gamma_M\Lambda^{p+1}(\ker\rho)\otimes \Xx^q(N)$$ where $\ud$ is the naive
differential. The result follows.
\end{proof}

The third sheet of the spectral sequence is trivially deduced from Proposition~\ref{p:E2}.
\begin{cor}\label{c:E3} Let $E$ be a Courant algebroid with split base. There
  is a canonical isomorphism  of bigraded algebras
 $E_3^{\bullet,\bullet} \cong E_2^{\bullet,\bullet}$.
\end{cor}
\begin{proof}
 Since $\Xx^q(N)$ is concentrated in even degrees for $q$, so is $E_2^{p,q}$
 by Proposition~\ref{p:E2}. Therefore, the differential $\ud_2:E_2^{p,q}\to
 E_2^{p+2,q-1}$ is necessarily $0$. Hence $E_3^{\bullet,\bullet} \cong
 E_2^{\bullet,\bullet}$.
\end{proof}

\bigskip

We  now prove the conjecture of Sti{\'e}non-Xu.  There is a canonical
morphism $\phi:H^\bullet_{naive}(E) \to H^\bullet_{std}(E)$, see~\cite{SX07}
and  Section~\ref{s:naive}.
\begin{cor}\label{c:trans} Let $E$ be a transitive Courant algebroid. Then
  the canonical map $\phi$ is an isomorphism $\phi:H^\bullet_{naive}(E)\cong 
 H^\bullet_{std}(E)$, {\em i.e.}, the Courant algebroid cohomology
coincides with the naive cohomology.
\end{cor}
\begin{proof}
A transitive Courant algebroid satisfies $D:=\img\rho=TM$.  Therefore, $E$ is
trivially with split base and
$M/D=\pt=N$.  Hence $\Xx^q(N)$ is non zero only for $q=0$, where it is
$\R$. Therefore, by Proposition~\ref{p:E2},  $E_2^{p,q}=0$ if $q\neq 0$. It
follows that all the higher differentials $\ud_{n\geq 3}$ are null. Thus the
cohomology of the Courant algebroid is isomorphic to $E_2^{\bullet,0}\cong
H^\bullet_{naive}(E)\otimes \R$.

\smallskip

Furthermore, by definition of the naive ideal $I$, the map $\phi$ preserves
the filtration by $I$. Thus $\phi$ passes to the spectral
sequence and coincides with the morphism of complexes
$\phi_1:(\Gamma(\Lambda^n \ker \rho),\ud)\cong  (E_1^{n,0},\ud_1)
\hookrightarrow (\oplus_{p+q=n} E_1^{p,q},\ud_1)$  on the first sheet of the
spectral sequence. By the first paragraph of this proof,  $\phi_1$ is a
quasi-isomorphism. Hence $\phi$ is indeed an isomorphism.
\end{proof}
\begin{rem}\label{rem:conjectureexact}
When $E$ is an exact Courant algebroid,  Corollary~\ref{c:trans} was obtained
by Roytenberg~\cite{Royt02}. Indeed, it was one of the examples motivating the
conjecture of Sti\'enon--Xu.
\end{rem}
\bigskip

For general Courant algebroids $E$ with split base, the spectral sequence
does not collapse on the sheet $E_2^{\bullet,\bullet}$ but is
controlled by a map from the vector fields on $N$ to the naive
cohomology of $E$ which we now describe.

It is a general fact from the spectral sequence theory that there exists a 
differential $\ud_3:E_3^{p,q}\to E_3^{p+3,q-2}$. In particular, $\ud_3$ 
induces an $E_3^{0,0}=\smooth(N)$-linear map 
\begin{equation} \label{eq:T3} 
  T_3:\X(N) \to H^3_{naive}(E) 
\end{equation} given as the composition 
$$ T_3:\X(N)\cong \Xx^2(N)\cong E_3^{0,2}\stackrel{\ud_3}\longrightarrow E_3^{3,
0}\cong H^3_{naive}(E) \;.$$ We call the map $T_3$ the \emph{transgression 
homomorphism} of the Courant algebroid $E$. Let $\X^{kil}(N)$ 
be the kernel of $T_3$ above (we 
like to think of elements of $\X^{kil}(N)$ as Killing vector fields 
preserving the structure function $H$). Note that $\X^{kil}(N)$ may be 
singular, {\em i.e.}, its rank could vary. We denote $\Xx^{kil, q}$ the 
space of ``symmetric Killing multivector fields'' $S^{q/2}_{\smooth(N)
}(\Xx^{kil}(N))$ with the convention that $\Xx^{kil, q}=\{0\}$ for odd 
$q$.

\begin{thm}\label{thm:split} The  cohomology of a Courant algebroid $E$ with
  split base is given by
\begin{equation}\label{eq:split}  H^n_{std}(E)
  \cong \bigoplus_{p+q=n} H^p_{naive}(E)/(T_3) \otimes \Xx^{kil, q} \;,
\end{equation} where $(T_3)$ is the ideal in $H^\bullet_{naive}(E)$ which
is generated by the image $T_3(\X(N))$ of $T_3$.
\end{thm}
\begin{proof} According to Proposition~\ref{p:E2} and Corollary~\ref{c:E3}, 
the third sheet of the spectral sequence is given by $E_3^{p,q}\cong 
H_{naive}^p(E)\otimes \Xx^q(N)$. Note that $\Xx^q(N)= \Gamma_N(S^{q/2}(TN))$ 
is generated as an algebra by its degree 2 elements. Since the differential 
$\ud_3:E_3^{p,q}\to E_3^{p+3,q-2}$ is a derivation, it is necessarily the 
unique derivation extending its restriction $T_3=\ud_3: \X(N)\to H_{naive}^3(E)
$. Since $\Xx^\bullet(N)=\Gamma_N(S^{q/2}(TN))$ is a free graded commutative 
algebra, the cohomology $E_4^{\bullet,\bullet}=H^\bullet(E_3^{\bullet,\bullet},
 \ud_3)$ is given by
\begin{eqnarray*}
  E_4^{\bullet,q} &\cong& H_{naive}^\bullet(E)/({\rm im}(T_3)) \otimes
    S^{q/2} (\ker T_3) \\
  &\cong & H^\bullet_{naive}(E)/(T_3) \otimes \Xx^{kil, q} \;.
\end{eqnarray*}
Now, it is sufficient to prove that all higher differentials $\ud_{r\geq 4}$ 
vanish. Since $\ud_4:E_4^{p,q}\to E_4^{p+4,q-3}$ is a derivation, it is 
completely determined by its restriction to the generators of $E_4^{\bullet,
\bullet}$ which lie in $E_4^{\bullet,0}$ and $E_4^{0,2}$. Thus, for obvious 
degree reasons, $\ud_4=0$ and similarly for all $\ud_{r>4}$. Therefore, the 
spectral sequence collapses at the fourth sheet\,: $E_\infty^{p,q}\cong E_4^{p,
 q}$.
\end{proof}
\begin{rem}\label{rem:algebra}
Theorem~\ref{thm:split} gives an isomorphism of vector spaces, not of
algebras in general. However,  by standard results on spectral sequences of
algebras, the  isomorphism~\eqref{eq:split} is an isomorphism of graded
algebras if the right hand side of~\eqref{eq:split} is free as a graded
commutative algebra.
\end{rem}
\begin{rem} Theorem~\ref{thm:split} implies that all the cohomological 
information of a Courant algebroid with split base is encoded in the 
transgression homomorphism $T_3$ together with the naive cohomology (and image 
of the anchor map). We like to think of $T_3$ as a family of closed 3-sections 
of $\ker \rho$ obtained by transgression from $\X(N)$. This idea is made more 
explicit in the case of generalized exact Courant algebroids in 
Section~\ref{s:gexsplit}. In that case, the transgression homomorphism is 
closely related to a generalization of the characteristic class of the Courant 
algebroids as defined by  {\v S}evera~\cite{SevLett}, that is the 
cohomology class of the structure 3-form (see Proposition~\ref{P:Severaclass} 
and Proposition~\ref{P:D=C}) parameterizing such Courant algebroids.
\end{rem}

\section{Generalized exact Courant algebroids with split base}\label{s:gexsplit}
In this section, we consider a generalization of exact Courant
algebroids. These Courant algebroids are parametrized by the cohomology class
of closed $3$-forms from which an explicit formula for the transgression
homomorphism $T_3$ can be given.

An \emph{exact}
Courant algebroid $E\to M$ is a Courant algebroid  such  that the following
sequence 
$$ 0\to T^*M \xrightarrow{\rho^*} E\xrightarrow{\rho} TM\to 0 $$
is exact.

Assume that $D:=\im\rho\subset TM$ is a subbundle, \emph{i.e.}, $E$ is
regular. Then the anchor maps surjectively $E\stackrel{\rho}\to D$ and its
dual $\rho^T:T^*M\to E^*\cong E$ factors through an injective map
$D^*\stackrel{\rho^*}\to E$ (again $E$ and $E^*$ are identified by the
pseudo-metric). 
\begin{vdef}\label{d:genexact} A regular Courant algebroid
  $(E,\rho,\ldots)$ is generalized exact if the following sequence of bundle
  morphisms over $M$
\begin{eqnarray}\label{eq:genexact} 0\to D^*\xrightarrow{\rho^*} E\xrightarrow{\rho} D\to 0 \end{eqnarray}
is exact.
\end{vdef}
By the above discussion, the only condition to check for the
sequence~\eqref{eq:genexact} to be exact is exactness in $E$.

\medskip

There is a simple geometric classification of exact Courant algebroids due to
{\v S}evera~\cite{SevLett} which extends easily to generalized exact Courant
algebroids as follows. Note that $D$ is Lie subalgebroid of $TM$. Given a Lie algebroid $D$, we write $(\Omega^\bullet(D),\ud_D)$ the complex of $D$-forms, $\Hal^\bullet(D)$ its cohomology and $Z^\bullet(D)=\ker(\ud_D)$ the closed $D$-forms~\cite{MacK87, CaWe99}.

Let $(E\to M,\rho,[.,.],\langle.,.\rangle)$ be a generalized exact Courant
algebroid  and assume given a splitting of the exact
sequence~\eqref{eq:genexact} as a pseudo-Euclidean vector bundle, that is, an
isotropic (with respect to the pseudo-metric) section $\sigma:D\to E$ of
$\rho$. Thus the section
$\sigma$ identifies $E$ with $D^*\oplus D$ endowed with its standard
pseudo-metric: $\langle \alpha\oplus X, \beta\oplus Y\rangle=
 \alpha(Y)+\beta(X)$. Since $\rho$ preserves the bracket, for any 
$X,Y\in \Gamma(D)$, one has  $$[\sigma(X),\sigma(Y)]=\sigma([X,Y]_{TM}) \oplus
\tilde{C}_\sigma(X,Y)$$ where $\tilde{C}_\sigma(X,Y) \in \rho^*(D^*)$. Let
$C_\sigma$ be the dual of $\tilde{C}_\sigma$, that is, for $X,Y,Z\in
\Gamma(D)$, we 
define $C_\sigma(X,Y,Z)=\langle \tilde{C}_\sigma(X,Y),Z\rangle$. It follows
from axiom~\eqref{ax:nSkew} and axiom~\eqref{ax:adInvar} of a Courant
algebroid (see Definition~\ref{def:Courant}) that $C_\sigma$ is
skew-symmetric. Moreover, by axiom~\eqref{ax:Leibn} and
axiom~\eqref{ax:adInvar}, $C_\sigma$ is $\smooth(M)$-linear. Thus $C_\sigma$
is indeed a $3$-form on the Lie algebroid $D$,\emph{ i.e.}, $C_\sigma\in
\Omega^3(D)$. Furthermore, the (specialized) Jacobi identity, \emph{i.e.},
axiom~\eqref{ax:Jacobi} implies that $C_\sigma$ is closed, that is,
$\ud_D(C_\sigma)=0$.
\begin{prop}[The {\v S}evera characteristic class] \label{P:Severaclass} Let 
$E$ be a generalized exact Courant algebroid.
\begin{enumerate}\item There is a splitting of the exact sequence~\eqref{eq:genexact} 
as a pseudo-Euclidean bundle; in particular there is an isotropic section 
$\sigma: D\to E$ of $\rho$.
\item \label{Severaclass} If $\sigma':D\to E$ is another isotropic section, 
then $C_\sigma -C_{\sigma'}$ is an exact $3$-form, that is, 
$C_\sigma -C_{\sigma'}\in \im \ud_D$.
\end{enumerate}
\end{prop}
\begin{proof}  The proof is the same as the one for exact Courant algebroids, 
for instance, see {\v S}evera Letter~\cite{SevLett} or~\cite[Section 3.8]{Royt99}.
\end{proof}
In particular the cohomology class $[C_\sigma]\in \Hal^3(D)$ is independent of 
$\sigma$. We will simply denote it $[C]$ henceforth.  We call the class $[C]\in \Hal^3(D)$ the 
{\v S}evera class of $(E\to M,\rho,[.,.],\langle.,.\rangle)$.

\smallskip

Given a closed $3$-form $C \in \Omega^3(D)$, one can define a bracket on the 
pseudo-Euclidean vector bundle $D^*\oplus D$ given, for $X,Y\in \Gamma(D)$,
$\alpha, \beta \in \Gamma(D^*)$, by the formula
\begin{equation}\label{eq:formbracket}
  [\alpha \oplus X, \beta \oplus Y] = \Lie_X\beta -\imath_Y\ud_D(\alpha)+C(X,Y,.) \oplus [X,Y]_{TM}.
\end{equation}
It is straightforward to check that this bracket makes the pseudo-Euclidean 
bundle $D^*\oplus D$ a Courant algebroid, where the anchor map is the 
projection $D^*\oplus D\to D$~\cite[Section 3.8]{Royt99}. Clearly its 
{\v S}evera class is $C$. Moreover two cohomologous closed 
$3$-forms $C,C'\in \Omega^3(D)$ yield isomorphic Courant 
algebroids~\cite[Section 3.8]{Royt99}. Therefore 
\begin{prop}[Analog of the {\v S}evera classification]\label{P:classification} 
Let $D$ be a Lie 
subalgebroid of a smooth manifold $M$. The isomorphism classes of generalized 
exact Courant algebroids with fixed image $\im \rho =D$ are in one to one 
correspondence with $\Hal^3(D)$, the third cohomology group of the Lie 
algebroid $D$.

The correspondence assigns to a Courant algebroid $E$ its {\v S}evera 
characteristic class $[C]$.
\end{prop}

\notneeded{\begin{rem}\PX{Really not needed}
By Proposition~\ref{P:Severaclass} and Proposition~\ref{P:classification}, 
generalized exact Courant algebroids are a special case of the (Courant 
algebroid associated to) quasi-Lie bialgebroids as defined by 
Roytenberg~\cite{Royt99}.
\end{rem}}

\begin{rem}\label{rem:Hgenexact} Generalized exact Courant algebroids are easy 
to describe in terms of the derived bracket construction. Indeed, since we can 
choose an isomorphism $E\cong D\oplus D^*$, the minimal symplectic realization 
of $E$ is isomorphic to $T^*[2]D[1]$. From the explicit 
formula~\eqref{eq:formbracket} for the Courant bracket, we found that the 
generating cubic Hamiltonian $H$ (encoding the Courant algebroid structure) is 
given, in our adapted coordinates, by 
\begin{equation} \label{eq:Hgenexact} 
  H=p_I\xi^I+\frac16C_{IJK}\xi^I\xi^J\xi^K 
\end{equation} where $C_{IJK}$ are 
the components of the {\v S}evera $3$-form $C$ induced by the splitting $\cong 
D\oplus D^*$.
\end{rem}

\medskip

Now assume that $E$ is a generalized exact Courant algebroid with split base 
$D\cong TL\times N$. Then, there is an isomorphism $\Omega^\bullet(D)\cong 
\Omega^\bullet(L)\otimes_\R\smooth(N)$ and the de Rham differential $\ud_D$ of 
the Lie algebroid $D$ is identified with $\ud_L \otimes 1$, the de Rham 
differential of the smooth manifold $L$, see Remark~\ref{rem:deRham}. \notneeded{In other 
words, the de Rham complex of the Lie algebroid $D$ is the de Rham complex of 
the manifold $L$ tensored by smooth functions on $N$.} In particular, $Z^3(D)
\cong Z^3(L)\otimes_\R\smooth(N) $. Furthermore, since $\ker \rho \cong 
D^*\cong T^* L\times N$, the naive complex $\big(\Gamma_M(\Lambda^\bullet \ker 
\rho),\ud\big)$ is isomorphic to the complex $\big(\Omega^\bullet(L)\otimes 
\smooth(N), \ud_L\otimes 1\big)$ of de Rham forms on $L$ tensored by smooth 
functions on $N$. Therefore $H_{naive}^\bullet(E)\cong H^\bullet_{DR}(L)\otimes 
\smooth(N)$, where $H^\bullet_{DR}(L)$ is the de Rham cohomology of the
manifold $L$.

\begin{ex} \label{e:genexactsplitbase}
Let $L$ and $N$ be two smooth manifolds and define $M:=L\times N$, 
$D:=TL\times N$ which is a subbundle of $TM$. Pick $\omega\in Z^3(L)$ any 
closed $3$-form on $L$ and $f\in \smooth(N)$ be any function on $N$. Then 
$C:=\omega \otimes f$ is a closed $3$-form on $D$. The $3$-form $C$ induces a 
generalized exact Courant algebroid with split base structure on $E:= D\oplus 
D^*$ where the Courant bracket is given by formula~\eqref{eq:formbracket}, the 
pseudo-metric is the standard pairing between $D$ and $D^*$ and the anchor map 
the projection $E\to D$ on the first summand. By 
Proposition~\ref{P:classification}, any generalized exact Courant algebroid 
with split base is isomorphic to such a Courant algebroid.
\end{ex}

\bigskip

\notneeded{\PX{You should write this paragraph much shorter.  No need to use such fancy
  language.}For a generalized exact Courant algebroid with split base $E$,
the transgression homomorphism $T_3:\X(N)\to H^3_{naive}(E)$  is easily
described in terms of the {\v S}evera class of $E$.}
 For any vector field $X\in \X(N)$ on $N$, there is the map 
\begin{equation} \label{eq:T3gen} \Omega_{naive}^\bullet(E)\cong
  \Omega^\bullet(L)\otimes \smooth(N) \xrightarrow{\; 1\otimes X\; }
  \Omega^\bullet(L)\otimes \smooth(N)\cong
  \Omega_{naive}^\bullet(E)
\end{equation} 
defined, for $\omega \in \Omega^\bullet(L)$ and $f\in \smooth(N)$ by
$(1\otimes X)(\omega\otimes f)=\omega \otimes X[f]$. 
Applying this map to the {\v S}evera $3$-form $C\in Z^3(L)\otimes
\smooth(N)\cong Z^3(D)$ yields the map  $$\X(N) \ni X\mapsto [(1 \otimes X)(C)]\in
H^3_{naive}(E)$$ which is well defined and depends only on the {\v S}evera
class of $E$ (and not on the particular choice of a $3$-form
representing it) by  Proposition~\ref{P:Severaclass}.
\begin{prop}\label{P:D=C}
 Let $E$ be a generalized exact Courant algebroid with {\v S}evera
 class $[C]\in \Hal^3(D)\cong H^3_{DR}(L)\otimes \smooth(N)$. The
 transgression homomorphism $T_3:\X(N)\to H^3_{naive}(E)$ is the map given,
 for any vector field $X\in \X(N)$, by 
$$T_3(X) = [(1 \otimes X)(C)] \in H^3_{naive}(E).$$
\end{prop}
\begin{proof} Fix a {\v S}evera $3$-form $C$ representing the {\v S}evera
  characteristic class. The transgression homomorphism is induced by the
  differential $\ud_3:E_3^{\bullet,\bullet}\to E_3^{\bullet+3,\bullet-2}$.
  Corollary~\ref{c:E3} and Proposition~\ref{p:E2} together with above remarks
  about generalized exact Courant algebroids give\,:
  $E_3^{p,q}\cong\Gamma(\Lambda^p\im\rho^*)\otimes\Xx^q(N).$  Therefore our
  adapted coordinates \ref{rem:lcoo} still apply and give on a local chart
  the map $\ud_3$ as
\begin{align*}
 \ud_3
  =& \frac16C_{IJK,L'}\xi^I\xi^J\xi^K\pfrac{}{p_{L'}}
\end{align*} \MG{Note that there are no further $\xi^A$, because of
  generalized exactness.}
By formula~\eqref{eq:Hgenexact}, the functions
$\frac16C_{IJK}\xi^I\xi^J\xi^K$ are given by the components of the
{\v S}evera $3$-form $C$ (identified with a function on 
$E[1]$ via $\langle.,.\rangle$ and pulled back to $\E$). Now the result
follows since $T_3$ is the restriction of $\ud_3$ to $E_3^{0,2}\cong \X(N)$.
\end{proof}

\medskip

 We denote $\Ann(C)$ the kernel of $T_3$, that is the vector fields $X$ on $N$ such that $[(1\otimes X)(C)]=0\in H^3_{naive}(E)$.  We also denote $\big((1\otimes \X(N))(C)\big)$ the ideal in $H^\bullet_{naive}(E)$ generated by the vector subspace $\im T_3 =\{(1\otimes X)(C) \, , \, X\in \X(N)\}$.
\begin{cor}\label{c:genexact} The Courant algebroid cohomology of a generalized exact Courant
 algebroid $E$ with split base $D\cong TL\times N$  is given by
\begin{align*} H^n_{std}(E)
  &\cong \bigoplus_{p+2q=n} H^p_{DR}(L)\otimes \smooth(N)/\big((1\otimes \X(N))(C)\big) \otimes_{\smooth(N)}
    S^{q}(\Ann(C))
\end{align*}
where $[C]$ is the  {\v S}evera class of $E$.
\end{cor}
\begin{proof} This is an immediate consequence of Theorem~\ref{thm:split}, 
Proposition~\ref{P:D=C} and the isomorphism $H^\bullet_{naive}(E)\cong 
H^\bullet_{DR}(L)\otimes \smooth(N)$.
\end{proof}

\begin{rem}\label{e:D3triv}
 Let $E$ be a generalized exact Courant algebroid with split base $D\cong 
 TL\otimes N$. Assume that the {\v S}evera class of $E$ can be represented by a 
 $3$-form $C\in \Omega^3(L)\otimes \gsmooth(N)$ which is constant as a function 
 of $N$, \emph{i.e.}, $C\in \Omega^3(L)\otimes \R \subset 
 \Omega^3(L)\otimes \gsmooth(N)$. Then, by Proposition~\ref{P:D=C} $T_3=0$ and 
 thus, by Corollary~\ref{c:genexact}, the cohomology of $E$ is 
 $$H^n_{std}(E) \cong \bigoplus_{p+q=n} H^p_{DR}(L)\otimes \Xx^{q}(N) \,.$$ 
 An example of such a 
 Courant algebroid is obtained as in Example~\ref{e:genexactsplitbase} by 
 taking $C=\omega \otimes 1$ for the {\v S}evera $3$-form, where $\omega$ is 
 any closed $3$-form on $L$.
\end{rem}

\begin{ex}  Let $G$ be a Lie group with a bi-invariant metric
$\langle.,.\rangle$. Then $G$ has a canonical closed $3$-form which  is the
Cartan $3$-form $C=\langle [\theta^L,\theta^L],\theta^L\rangle$ where
$\theta^L$ is the left-invariant Maurer-Cartan $1$-form\MG{ and bracket and 
metric are understood on the Lie algebra component}. Note that, by
ad-invariance of  $\langle.,.\rangle$, $C$ is bi-invariant (hence closed).
\MG{That's also what Roytenberg said.  Seems I'm the only one who doesn't know 
this triviality.}

Thus, by Example~\ref{e:genexactsplitbase}, there is a generalized exact
Courant algebroid structure on $G\times N$ with {\v S}evera class $[C\otimes
  f]$ for any manifold $N$ and function $f\in\smooth(N)$. This example (for
$N=\{*\}$) was suggested by  Alekseev (see also \cite[example
  5.5]{Royt02}).  Explicitly, the Courant algebroid is $E=(\g\oplus\g)
\times N\to G\times N$. The structure maps are given by\,:
\begin{align*}
  \langle X\oplus Y, X'\oplus Y'\rangle = \langle X,X'\rangle -\langle
  Y,Y'\rangle  \\
  \rho: E\to TG\boxplus TN: (X\oplus Y,g,n) \mapsto (X^l-Y^r)(g)\boxplus0
\intertext{where the superscript $^l$ (resp.\ $^r$) means that an element of
  the Lie algebra $\g\cong T_eG$ is extended as a left (right) invariant
  vector field on $G$. The bracket is given, for $(g,n)\in G\times N$ and 
  $X,X',Y,Y'\in\g\subset\Gamma(\g\times G),$ by
} 
  [X\oplus X', Y\oplus Y']_{(g,n)} = ([X,X']\oplus  f(n)[Y,Y'])_{(g,n)}&
\end{align*}  Choosing the splitting $\sigma(Z,g,n):=(Z\oplus-\Ad_gZ,g,n)$, one 
finds that the Cartan $3$-form $C$ is indeed the {\v S}evera 
class of $E$.

Now assume $G$ is a compact simple Lie group, then $C$ spans $H^3(G)$. If we 
take $N=\R$ and $f$ to be constant, then, by 
Remark~\ref{e:D3triv}, the cohomology of $E$ is $H^\bullet_{std}(E)\cong H^\bullet(G)
\otimes \X^\bullet(\R)$, thus two copies (in different degrees) of the 
cohomology of $G$. Now take $f\in \smooth(\R)$ to be $f(t)=t$. Then $\Ann(C)$ 
is trivial and $ H^\bullet_{std}(E)= H^\bullet(G)/(C)$. Note that $H^\bullet(G)$
is the exterior algebra $H^\bullet(G)\cong 
 \Lambda^\bullet(C,x_2,\dots,x_r)$, thus $H^\bullet_{std}(E)\cong 
\Lambda^\bullet(x_2,\dots,x_r)$ as an algebra, see Remark~\ref{rem:algebra}. 
\end{ex}

\appendix
\section{Graded geometry} \label{app:graded}
In this appendix we recall some basics of super and graded geometry. For a good and detailed introduction to supermanifolds see, for instance,~\cite{Vor91}.

\def\Csuper{{\ensuremath{\C^\infty}}}

\begin{vdef}A \emph{supermanifold} $\E$ of dimension $p|q$ is a
  smooth manifold $M$ of dimension $p$ together with a sheaf, denoted $\Csuper(\E)$, of $\Z_2$-graded
  $\C^\infty(M)$-algebras locally of the form
  $\C^\infty(U)\otimes_\R\Lambda^\bullet(\R^q)$.
\end{vdef}
$\smooth(\E)$ is called the sheaf of super functions on $\E$.
\begin{ex} \label{ex:ovb} A standard example of a (non-trivial) super
  manifold is the odd vector bundle $\Pi E$ associated to a smooth vector
  bundle $E\to M$.  The supermanifold  $\Pi E$  is the  manifold $M$
  together with the structure sheaf 
$\Gamma_M(\Lambda^\bullet E^*)$. Locally, $\smooth(\E)$ is  generated by $C^\infty(U)$ and a local frame
$\xi^1,\ldots,\xi^k$ of $E^*$, where $k$ is the rank of $E$.
\end{ex}
Example~\ref{ex:ovb} is fundamental because of the well known
\begin{thm}[Batchelor] Every supermanifold can be realized as an odd vector
bundle.
\end{thm}
\PX{Give a reference for the thm/proof\,!}\MGI{done.}For a proof of this see, 
e.g. \cite{Vor91}.  In particular,  the functor $\Pi:Vect\to SMan$ (given by
Example~\ref{ex:ovb}) from the category of vector bundles to the category of
supermanifolds is surjective on the objects. However these two categories are
not equivalent, since the category of supermanifolds has more morphisms than
the category of vector bundles~\cite{Vor91}.

Vector bundles on supermanifolds are defined analogously to smooth vector
bundles, using  $\Z_2$-graded vector spaces instead of mere non graded vector
spaces. In particular, the tangent bundle of a supermanifold $\E$ is the
space of graded derivations of the structure sheaf $\smooth(\E)$.

\subsection{Graded manifolds}
For graded manifolds we refer to \cite[section 4]{Vor01}, \cite[section
  2]{Sev01b} or \cite[section 2]{Royt02}.

The difference between a supermanifold and a graded manifold lies in an
additional data for a graded manifold -- the \emph{Euler vector field} -- as
well as a cover of compatible $\Z$-graded charts.  Before
recalling the formal definition, let us first consider a fundamental example
along the lines of Example~\ref{ex:ovb}.

\begin{ex} \label{ex:evf} Let $E\to M$ be a smooth vector bundle and
  consider the supermanifold $\E=\Pi T \Pi E$. An Euler vector field
  $\epsilon$ can be defined, in coordinates induced by charts of $E\to M$, by
  the formula 
$$ \epsilon= 2\theta^a\pfrac{}{\theta^a}\,
  +v^i\pfrac{}{v^i}+\xi^a\pfrac{}{\xi^a}\,
$$
where $x^i, \xi^a$ are coordinates on $M$ and the fibers of $E$ and $v^i,
\theta^a$ are their correspondents on $TE$. Note that $\epsilon$ is well 
defined because for these adapted charts we have the following transition 
functions\,:
\begin{align*}  \widetilde{x^i} &= \widetilde{x^i}(x)  \\
  \widetilde{\xi^a} &= M^{\tilde a}_b(x)\xi^b \\
  \widetilde{v^i} &= \pfrac{\widetilde{x^i}}{x^j}v^j \\
  \widetilde{\theta^a} &=M^{\tilde a}_b\theta^b +M^{\tilde a}_{b,i}\xi^b v^i.
\end{align*}
Now we restrict $\smooth(\E)$ to functions which are polynomial in $\theta^a$s.
A crucial observation is that the adapted coordinates
are eigenfunctions of this Euler vector field $\epsilon$, namely
$\epsilon\cdot  x^i=0, \epsilon \cdot \xi^a= \xi^a, \epsilon \cdot  v^i=v^i$ and
$\epsilon\cdot \theta^a=2\theta^a$.  Moreover the structure sheaf
$\smooth(\E)$ is $\Z$-graded and locally (over such a chart) generated by these
coordinates (as smooth functions) in all coordinates of degree 0, the free
exterior algebra in all odd coordinates, and a polynomial
algebra in all even coordinates not of degree 0.
Therefore, the eigenvalues of $\epsilon$ are called degrees, and are
compatible with the $\Z_2$-grading of $\smooth(\E)$, \emph{i.e.},
eigenfunctions for odd eigenvalues are odd functions in the $\Z_2$-grading
and similarly for even eigenvalues.
The supermanifold $\E$ together with this Euler vector field and the charts
induced by $E\to M$ is an example of a graded manifold, denoted $T[1]E[1]$.
\end{ex}
\deprecated{The properties of the eigenfunctions of the Euler vector field in
Example~\ref{ex:evf} above are the ones defining a graded manifold\,:}

\begin{vdef}[(integer) graded manifold]\label{d:Nmfd} Let $\E$ be a fixed
  supermanifold with a fixed even vector field $\epsilon$.
\begin{enumerate}\item A chart of $\E$ is called $\Z$-graded iff its
  coordinates are eigenfunctions of $\epsilon$ with integer eigenvalues.

  The \emph{structure sheaf} $\gsmooth$ of $\Z$-graded functions over this
  chart is the $\Z$-graded algebra generated by these coordinates,
  i.e.\ smooth functions in the 
  coordinates of degree 0, the free exterior algebra in the odd coordinates,
  and the free algebra of symmetric polynomials in the even coordinates not
  of degree 0, with the $\Z$-grading given by $\epsilon$.
\item A $\Z$-graded atlas of $\E$ is an open cover with $\Z$-graded
  charts such that the transition functions between them preserve
  the $\Z$-grading and are constituted of $\Z$-graded functions.  Especially the
  number of coordinates in each degree is the same on all charts.

\item An (integer) graded manifold $\E$ is a supermanifold $\E$ together with
  a fixed vector field $\epsilon$, called the Euler vector field, and a
  (maximal) $\Z$-graded atlas.
\end{enumerate}
\end{vdef}
The graded manifold is said to be \emph{non negatively graded}
  if all coordinates are of non negative integer degree.

\MG{Maybe this follows now trivially from the ``generating'' rules.} Even
though the $\Z_2$-grading (even/ odd) is not required to be related
to the integer grading (by the Euler vector field) in the definition, in many
practical examples, the two gradings are compatible, \emph{i.e.}, odd
functions with respect to the
integer grading are exactly the odd with respect to the $\Z_2$-grading.  All
graded manifolds in this paper fulfill this restriction.
Non negatively graded manifolds for which the two gradings are compatible are
also called  N-manifolds~\cite{Royt02}.


\begin{ex}
Given a smooth vector bundle $E\to M$ and an integer $n\neq 0$, we can form
the graded manifold $E[n]$ with the canonical Euler vector field which
assigns the degree $n$ to fiber-linear functions as in
Example~\ref{ex:evf}.
\end{ex}
There is an obvious notion of vector bundles over
graded manifolds similar to the notion of vector bundles over  super
manifolds.  If $\E\to \mathcal{M}$ is  such a vector bundle, one denotes by
$\E[n]$, the graded vector bundle obtained by shifting the  fiber degrees
by $n$.

For instance, any graded manifold $\E$ has a tangent bundle $T\E$
defined as the graded vector bundle of (graded) derivations of
$\gsmooth(\E)$.
Similarly, its cotangent bundle $T^*\E$ is its
$\gsmooth(\E)$-linear dual.

A sequence of graded manifolds $ \C\to \E\to\mathcal{F}$ is  said to be \emph{exact} if the sequence of sheaves
$ \gsmooth(\C)\leftarrow\gsmooth(\E)\leftarrow\gsmooth(\mathcal{F})$ is
exact.   This provides a bridge between  graded geometry and homological
algebra.  \notneeded{Analog to $\smooth$ in the smooth category $\gsmooth$ is a
contravariant functor, i.e.\ maps between graded manifolds induce (by
definition) maps between the structure sheaves in the opposite direction.}

\subsection{Graded Poisson and symplectic manifolds}
There is a de Rham differential
$\ud_\E:\gsmooth(\E)\to \Omega^1(\E):=\Gamma_\E(T^*[1]\E)$ for graded
manifolds~\cite{Vor91, Royt02}, generalizing the de Rham differential for
smooth (non graded) manifolds,
which uniquely extends as a derivation to $\Omega^\bullet(\E):=
\Gamma_\E(S^\bullet (T^*[1]\E))$.
\notneeded{ to
 the graded forms, i.e.\ sections into $\Lambda^\bullet T^*\E$.  Strictly
 speaking this is a double graded vector bundle, but the most important
 grading is its total grading which corresponds to the notation $S^\bullet
  T^*[1]\E$, where again $S^\bullet$ means the graded symmetric powers.}

\begin{vdef}  A \emph{graded Poisson
  manifold} of degree $d$, $(\E,\{.,.\})$ is a graded manifold $\E$ together
  with a  bracket of degree $d$ on its structure sheaf $\gsmooth(\E)$
  satisfying, for all homogeneous functions $f,g,h\in \gsmooth(\E)$,
\notneeded{, i.e.\ there is a fixed integer $d$ (called the degree of
the bracket) such that with the notation $A^\bullet:=\gsmooth(\E)$, $\{.,.\}:
A^p\otimes_\R A^q\to A^{p+q+d}$ and for all homogeneous functions $f,g,h\in
A^\bullet$:}
\begin{align}
  \{g,f\} &= (-1)^{(|f|+d)(|g|+d)+1} \{f,g\}  \label{eq:gSkew}\\
  \{f,gh\} &= \{f,g\}h +(-1)^{(|f|+d)|g|}g\{f,h\}  \label{eq:gLeibn}\\
  \{f,\{g,h\}\} &= \{\{f,g\},h\} +(-1)^{(|f|+d)(|g|+d)}\{g,\{f,h\}\}  \label{eq:gJacobi}
\end{align}

  A graded \emph{presymplectic manifold} $(\E,\omega)$ is a graded manifold
  with a homogeneous $\ud_\E$-closed 2-form $\omega$.
 
  A \emph{graded symplectic manifold} $(\E,\omega)$ is a graded presymplectic
  manifold with $\omega$ of total degree $d+2$ such that the map of vector
  bundles $\omega^\flat:T\E\to T^*\E$  is
  non-degenerate. It is in particular a graded Poisson manifold of degree
  $-d$.\MG{$\omega$ of degree$>2$ means that $\{.,.\}$ reduces the
    degree.  Of course $\deg\omega=2$ would mean the ungraded case.}
\end{vdef}

Graded Poisson manifolds  occur naturally  in field-theory, more particularly
in the context of the AKSZ or BV-formalism.  \notneeded{They are also the
  (graded) smooth analog of graded Poisson algebras.}  As in the non-graded
case, the Poisson bracket is  encoded by a homogeneous Poisson bivector
$\Pi$, \emph{i.e.}, living in $\Gamma_\E(S^2(T[-1]\E))$, satisfying 
$[\Pi,\Pi]_{Sc}=0$, where $[.,.]_{Sc}$ is the  graded analogue of the
Schouten bracket.
\notneeded{ an extension of the Lie bracket on the tangent space.  Its parity
  depends on the total grading of the graded skewsymmetric multivector
  fields.

  The meaning of the second ones lies in generalizing the third ones and in
  particular in the fact that for a graded map \phi:$\mathcal{F}\to\E$ into a
  (pre)symplectic manifold $(\E,\omega)$
  the pullback of the 2-form is still a closed 2-form (pullback a morphism of
  graded de Rham complexes), but in general no longer non-degenerate.

  The importance of the third ones is that they are at the same time a
  special case of the first ones, which has simpler properties.  (Analog to
  the ungraded case given a non-degenerate homogeneous 2-form $\omega$ then
  it is closed under the de Rham differential iff its inverse, a homogeneous
  bivector, is Poisson.)
}

\begin{ex}\begin{enumerate}
\item Let $(E, \langle .,.\rangle)$ be a pseudo-Euclidean  vector bundle over
  $M$. Then $E[1]$ has a canonical degree $-2$ Poisson bracket defined on two
  sections $\alpha,\beta\in \Gamma_M(E[1])$ by $\{\alpha,\beta\}=\langle
  \alpha, \beta\rangle$. Note that if $f,g\in \smooth(M)\cong
  {\gsmooth(E[1])}^0$, then $\{f,\alpha\}$ and $\{f,g\}$ are necessarily
  zero for degree reasons.  Note that $(E[1], \{ .,.\})$ is not a graded
  symplectic manifold.

  Similarly, one defines a degree $-2n$ Poisson manifold structure on $E[n]$.
\item The (shifted) cotangent bundles $T^*[k]\E$ is a  graded symplectic manifold with symplectic form
   $\omega=\ud_{T^*[k]\E}\theta$ induced by the 
   symplectic potential/Liouville-form $\theta$ similar to the smooth case.
  Here $\omega$ has total degree $k+2$.

  Let $z^\alpha$ be coordinates of $\E$ of degree $|\alpha|$, 
  $p_\alpha$ their duals of degree $k-|\alpha|$ and $\ud:=\ud_{T^*[k]\E}$ the
  de Rham differential on $T^*[k]\E$.  The Liouville form reads locally
  as  $ \theta := p_\alpha \ud z^\alpha$ of homogeneous degree $k+1$ and
  the symplectic form reads locally $$ \omega = \ud p_\alpha \ud z^\alpha $$
  which is obviously closed (under $\ud$) and also non-degenerate.  The Poisson bracket on the coordinates reads as
$$  \{p_\alpha,z^\beta\} = \delta_\alpha^\beta\quad
  \{z^\alpha,z^\beta\} = 0 = \{p_\alpha,p_\beta\}
$$ and has degree $-k$.

\notneeded{\item \def\g{\mathfrak{g}}Another example of a graded Poisson
  manifold that is no
  graded symplectic manifold would be $\g^*[k]$ for an ungraded Lie algebra
  $\g$ and even $k$, or just $\g^*$ for a graded Lie algebra $\g$.
}
\end{enumerate}
\end{ex}
\notneeded{In terms of algebra $\gsmooth(\E)$ with such a Poisson bracket is called a
$d$-bracket algebra.  Note that the usual (ungraded) Poisson bracket is
included for $k=0$, especially even.  If you pick $k$ odd, then the bracket
$\{.,.\}$ is a Gerstenhaber-bracket.}

\notneeded{A standard result in graded geometry would be that the fiber product of
 graded presymplectic morphisms is a graded presymplectic manifold if it is a
 graded manifold.  In order to get a graded symplectic manifold you need to
 check non-degeneracy of the graded presymplectic structure on the product.

 Unfortunately the situation for Poisson maps is more intricate, because already
 in the smooth case the fiber product generally has no canonic Poisson
 structure.
}

{\v S}evera has given the following example of a symplectic manifold which 
allows the derived bracket construction of Courant algebroids, see 
Section~\ref{s:defn}.
\begin{prop} Let $(E\to M, \langle .,. \rangle)$ be a pseudo-Euclidean vector bundle.  The
  fiber product $E[1]\times_{(E\oplus E^*)[1]} T^*[2]E[1]$ is a graded
  symplectic manifold where the symplectic form has total degree $2+2$.
\end{prop}
\begin{proof} Let us denote $g=\langle .,. \rangle$ the pseudo-metric and 
$g^\sharp:E[1]\to E^*[1]$ the isomorphism it induces. Let $i:E[1]\to (E\oplus 
E^*)[1]$ be the map $\psi\mapsto \psi\oplus \frac12g^\sharp\psi$. It is easy 
to check that $i$ is a morphism of graded Poisson manifolds. Thus the pullback 
along $i:E[1]\to (E\oplus E^*)[1]$ of the symplectic form on $T^*[2]E[1]$ 
gives a structure of graded presymplectic manifold on $E[1]\times_{(E\oplus E^*)
[1]} T^*[2]E[1]$. Let $\omega$ be the closed $2$-form. Now you need 
to check that $\omega$ is non-degenerate. It follows since, in the adapted 
coordinates from section~\ref{s:defn}, $\omega$ and the associated bivector 
$\Pi$, have the following explicit expression\,:
\begin{align*}
  \omega &= \ud p_i\wedge \ud x^i +\frac12g_{ab}\ud\xi^a\ud\xi^b \\
  \Pi &= \pfrac{}{p_i}\wedge\pfrac{}{x^i}+\frac12g^{ab}\pfrac{}{\xi^a}\cdot\pfrac{}{\xi^b}
\end{align*}
\end{proof}


\end{document}